\documentclass{gtart}
\input gtmonout
\volumenumber{1}
\volumeyear{1998}
\volumename{The Epstein birthday schrift}
\pagenumbers{365}{381}
\published{27 October 1998}
\received{2 September 1997}
\papernumber{18}

\usepackage[rokicki,hideboxes]{boxedeps}
\usepackage{labelfig}

\def\Conway{\BoxedEPSF{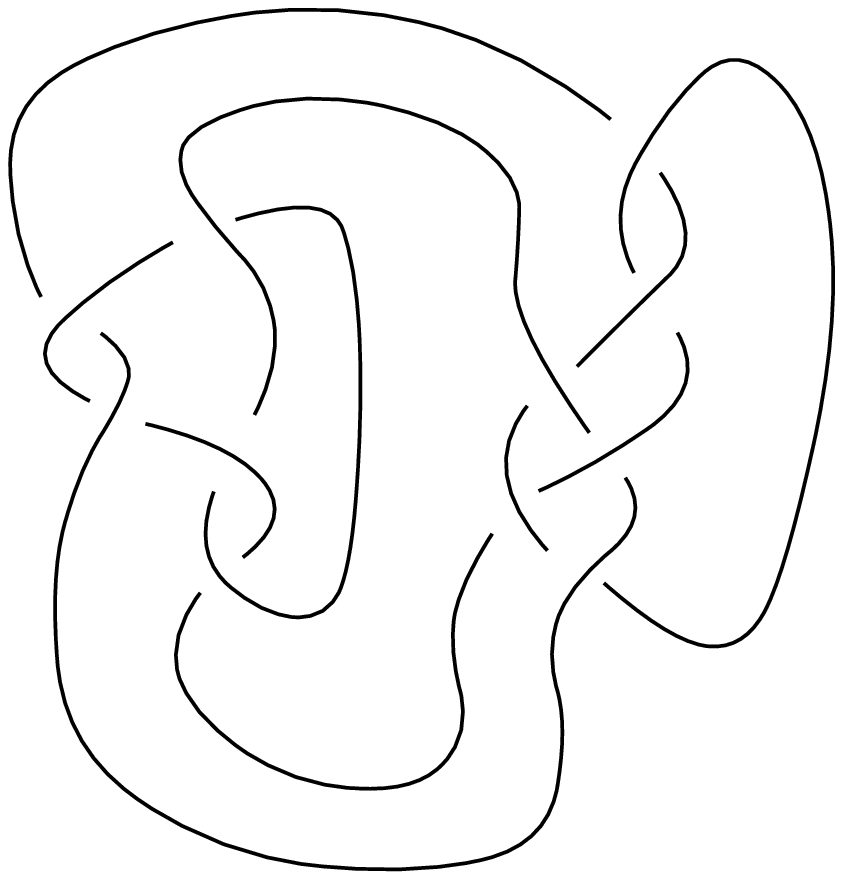 scaled  350}}
\def\KT{\BoxedEPSF{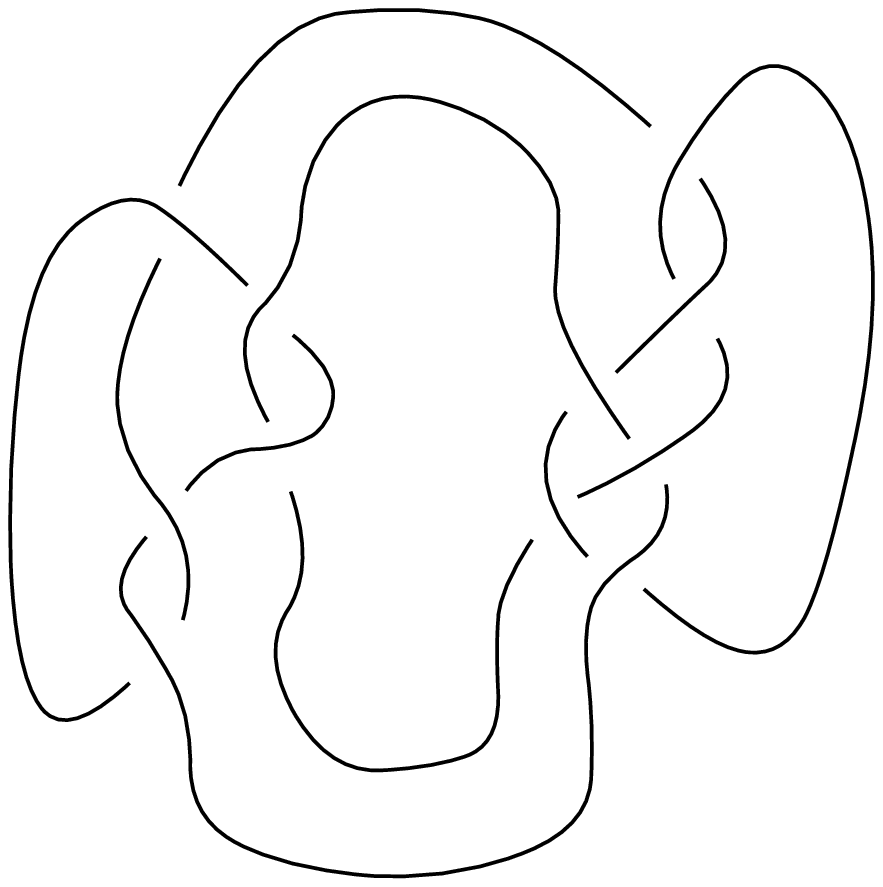 scaled  350}}

\def\Twoone{\BoxedEPSF{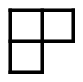 scaled  500}}
\def\Twotwo{\BoxedEPSF{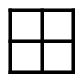 scaled  500}}
\def\Threeone{\BoxedEPSF{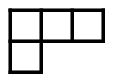 scaled  500}}
\def\Threetwo{\BoxedEPSF{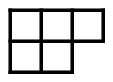 scaled  500}}
\def\Fourone{\BoxedEPSF{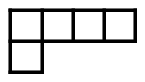 scaled  500}}

\def\Skeinelement{\BoxedEPSF{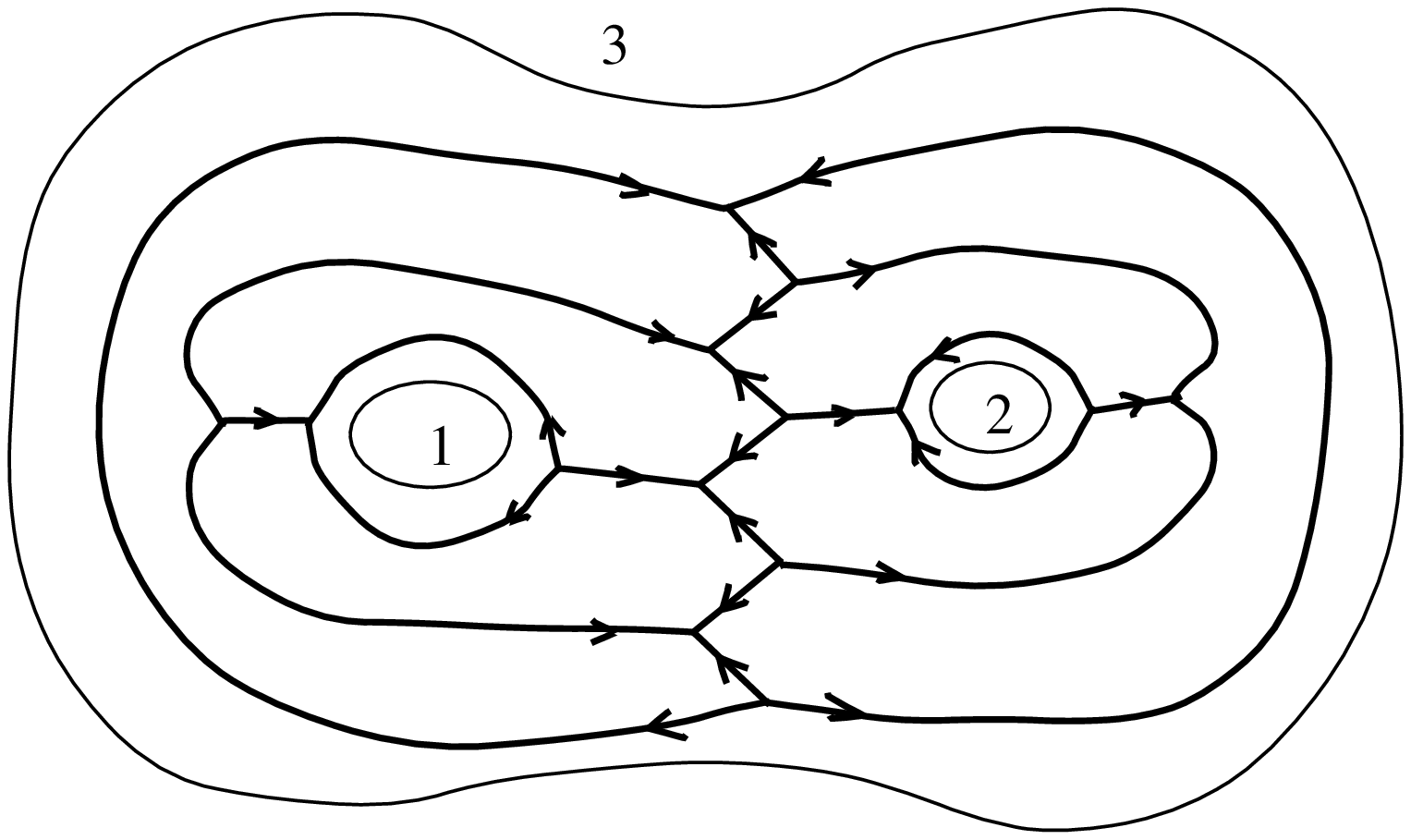 scaled  500}}
\def\Hexnet{\BoxedEPSF{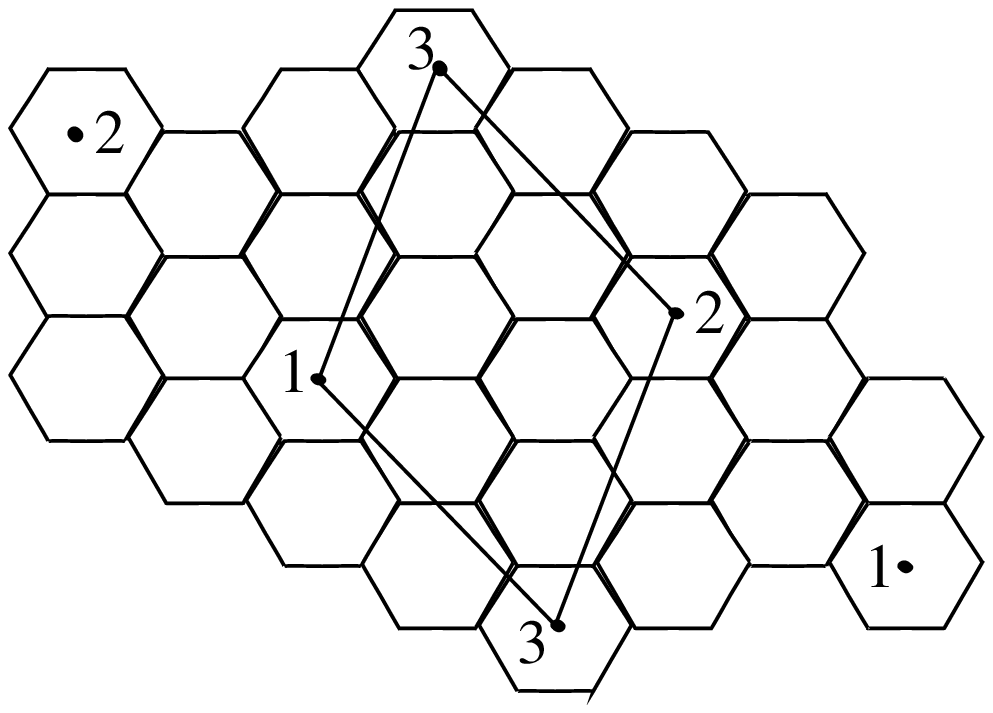 scaled  500}}

\def\Oplus{\bigoplus}
\def\ds{\displaystyle}

\newtheorem{thm}{Theorem}[section]

\begin{document}

\title{Mutants and $SU(3)_q$ invariants}         
%\shorttitle{Mutants and SU(3)_q invariants}
\asciititle{Mutants and SU(3)_q invariants}

\author{H\thinspace R Morton\\H\thinspace J Ryder} 
\shortauthors{Morton and Ryder}
\asciiauthors{H R Morton and H J Ryder}

\address{Department of Mathematical Sciences, University of Liverpool\\
Liverpool L69 3BX, England}

\email{h.r.morton@liv.ac.uk}  

\begin{abstract}
Details of quantum knot invariant  calculations using a specific
$SU(3)_q$--module are given which distinguish the Conway
and Kinoshita--Teresaka pair of mutant knots. Features of
Kuperberg's skein-theoretic techniques for $SU(3)_q$ invariants in the context
of mutant knots are also discussed.
\end{abstract}

\asciiabstract{%
Details of quantum knot invariant  calculations using a specific
SU(3)_q-module are given which distinguish the Conway
and Kinoshita-Teresaka pair of mutant knots. Features of
Kuperberg's skein-theoretic techniques for SU(3)_q invariants in the context
of mutant knots are also discussed.}



\maketitle

\section{Introduction} In previous studies of invariants derived from the
Homfly polynomial, or equivalently from the unitary quantum groups, it was
noted that no invariant given by a module over $SU(3)_q$  was known to
distinguish a mutant pair of knots. Indeed, any quantum group module whose
tensor square has no repeated summands determines a knot invariant which fails
to distinguish mutants \cite{MC}. A table of invariants which fail to
distinguish mutants was drawn up in \cite{MC}, using this and other evidence.
Direct Homfly polynomial calculations showed that a certain irreducible
$SU(N)_q$ invariant, coming from the module with Young diagram \Twoone, could
distinguish between some mutant pairs for $N \ge 4$, although not for $N=3$.
These calculations also exhibited a Vassiliev invariant of finite type 11 which
distinguishes some mutant pairs. The calculations left open the possibility
that $SU(3)_q$ invariants might never distinguish mutant pairs. 

In this paper we give details of calculations with a specific $SU(3)_q$--module
which result in different invariants for the Conway and Kinoshita--Teresaka pair
of mutant knots. We also consider some features of Kuperberg's skein-theoretic
techniques for $SU(3)_q$ invariants in the context of mutant knots.

Much of this work was carried out in 1994--95, while the second author was
supported by EPSRC grant GR/J72332.

\subsection{Background}

 The term {\it mutant\/} was coined by Conway, and  refers to the following
general construction.

Suppose that a knot $K$ can be decomposed into  two oriented  $2$--tangles $F$
and $G$ as shown in  figure 1. 

\begin{figure}[htbp]
\cl{%
\small
\SetLabels
\E(.05*.68) $K\ =$\\
\E(.18*.68) $F$\\
\E(.39*.68) $G$\\
\E(.6*.7) $K'\ =$\\
\E(.74*.7) $F'$\\
\E(.95*.7) $G$\\
\E(.05*.2) $F'\ =$\\
\E(.26*.2) $F$\\
\E(.42*.2) or\\
\E(.57*.2) $F$\\
\E(.75*.2) or\\
\E(.88*.2) $F$\\
\endSetLabels
%\ShowGrid
\AffixLabels{\BoxedEPSF{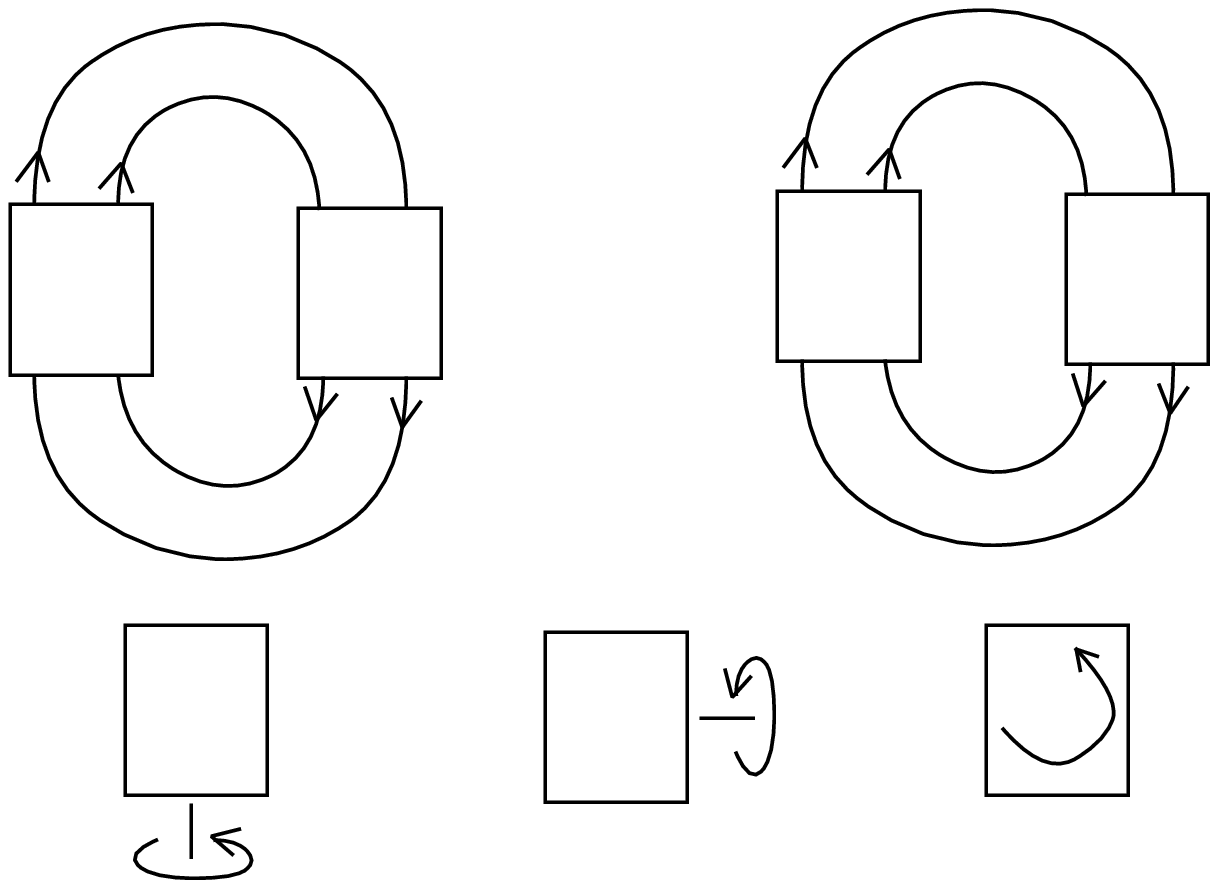 scaled  500}}}
\cl{\small Figure 1}
\end{figure}

 A new knot $K'$ can be formed by replacing the  tangle $F$ with the tangle
$F'$ given by rotating 
$F$ through $\pi $ in one of three ways,  reversing its string orientations if
necessary.  Any of these three knots $K'$ is called a  {\it mutant\/} of $K$. 

The two  $11$--crossing knots with trivial  Alexander polynomial found by 
Conway and  Kinoshita--Teresaka are the best-known example of  mutant knots.
They are shown  in figure 2.

\begin{figure}[htb]
\cl{\small $C\ =$\quad \Conway\ ,\qquad $KT\ =$ \quad  \KT\ .}
\smallskip 
\cl{\small Figure 2}
\end{figure}

 It is clear from figure 2  that the knots $C$ and $KT$ are mutants, and the
consituent tangles $F$ and $G$ are both given from a $3$--string braid by
closing off one of the strings.

The simplest $SU(3)_q$ invariant not previously known to agree on mutant pairs
is given by the 15--dimensional irreducible module with Young diagram \Threeone.
The Homfly polynomial of the 4--parallel with $z=s-s^{-1}$ and $v=s^N$ is a sum
of 4--cell invariants for $SU(N)_q$. When $N=3$ it is known that all 4--cell
invariants except that for \Threeone\ agree on mutants. Thus   the Homfly
polynomial of the 4--parallel, with the substitution $z=s-s^{-1}$ and $ v=s^3$,
agrees on mutants if and only if the $SU(3)_q$ invariant for \Threeone\ agrees
on mutants.

 Equally, the same substitution  in the Homfly polynomial of the satellite
consisting of the parallel with 3 strings, two oriented in one direction and
one in the reverse direction, gives the sum of certain 4--cell invariants for
$SU(3)_q$, because the dual of the fundamental module,  used to colour the
reverse string,  is given by using the Young diagram with a single column of
two cells. Then  the Homfly polynomial of the 3--parallel with one reverse
string, after the substitution $z=s-s^{-1},\; v=s^3$ agrees on mutants if and
only if  the $SU(3)_q$ invariant for \Threeone\ agrees on mutants.

Kuperberg's combinatorial methods for handling $SU(3)_q$ invariants seemed to
us for a while to offer a chance that the behaviour of $SU(3)_q$ would follow
that of $SU(2)_q$. We explored the $SU(3)_q$ skein of the pair of pants, based
on Kuperberg's combinatorial techniques, in the hope of proving this. An 
analysis of this skein is given later, as it has a geometrically appealing
basis, whose first lack of symmetry again pointed the finger at the reversed
3--parallel as the first potential candidate for distinguishing some mutant
pairs.

\subsection{Choice of calculational method}

 We did not pursue the Kuperberg skein calculations for these parallels of
Conway and Kinoshita--Teresaka. Although we contemplated briefly such an
approach it seemed difficult to use computational aids in dealing
with combinatorial skein diagrams once the number of crossings to be resolved
grew beyond easy blackboard calculations, as no computer implementation of the
graphical calculations in this skein was available to us. 

While we could, in principle, have calculated the Homfly polynomial of the
3--parallel of Conway's knot with one  string reversed there is
a considerable problem in computation of Homfly polynomials of links with
a large number of crossings. A number of computer programs will calculate the
Homfly polynomials of general links. Mostly these rely on implementation
of the skein relation, and the time required grows exponentially with
the number of crossings. Such programs include those by Ochiai, Millett
and Hoste. They will work up to the order of maybe 40 or even 50
crossings but  slow down rapidly after that. In the
application needed for this paper we have to deal with the
3--string parallel to Conway's knot with two strings in one direction
and one in the other, which gives a link with 99 crossings. Even if the
calculation is restricted to dealing with terms up to $z^{13}$ only, or some
similar bound, these programs are unlikely to make any impact on the
calculations.

There does exist a program, developed by Morton and Short
\cite{MS}, which can handle links with a large numbers of
crossings, under some circumstances. This is   based on the Hecke algebras, but
 it requires a braid presentation of the link on a restricted number
of strings; in practice 9 strings is a working limit, although in favourable
circumstances it can be enough to break the link 
 into pieces which meet this  bound more locally.   Unfortunately the reverse orientation of
one string which is needed in the present case means that any braid
presentation for the resulting link  falls well outside the
limitations of this program.

In \cite{MC} the Hecke algebra calculations on 3--string parallels with all
strings in the same direction could be carried out in terms of  9--string
braids, and  lent themselves well to an effective truncation to restrict the
degree of Vassiliev invariants which had to be calculated. The alternative 
possibility here  of using 4 parallel strings, all  with the same orientation,
faces the uncomfortable growth of these calculations from 9--string to 12--string
braids, entailing a growth in storage from 9! to 12! for a calculation which was
already nearing its limit. There are also
 almost twice as many crossings ($11\times16$), as well as a similar factorial
growth in overheads for the calculations.

We consequently did not pursue Homfly calculations any further.
Instead we returned
to the
$SU(3)_q$--module calculations and made explicit computations for the invariants
of the knots $C$ and $KT$ when coloured by the $15$--dimensional module
$V_{\Threeone}$, using the following scheme. This approach has the merit of
focussing directly on the key part of the $SU(3)_q$ specialisation, rather than
using the full Homfly polynomial on some parallel link. We give further details
of the method later. 

When each of the knots $C$ and $KT$ is coloured by the $SU(3)_q$--module $V_{\Threeone}$
the two constituent tangles $F$ and $G$ will be represented by an endomorphism
of the module $V_{\Threeone}\otimes V_{\Threeone}$. To calculate the invariant
of the knot, presented as the closure of the composite of the two 2--tangles, we
may compose the endomorphisms for the two 2--tangles, and then calculate the
invariant of the closure of the composite tangle in terms of the resulting
endomorphism. Let us suppose that
$V_{\Threeone}\otimes V_{\Threeone}$ decomposes as a sum $\Oplus
a_{\nu}V_{\nu}$ of irreducible modules, where $a_{\nu}\in {\bf N}$ and
$a_{\nu}V_{\nu}$ denotes the sum of all submodules which are isomorphic to
$V_{\nu}$. Any endomorphism then maps each isotypic piece $a_{\nu}V_{\nu}$ to
itself. It is convenient to regard each isotypic piece as a vector space of the
form $W_{\nu}\otimes V_{\nu}$, where $W_{\nu}$ has dimension $a_{\nu}$, and can
be explicitly identified with the space of highest weight vectors for the
irreducible module $V_{\nu}$ in $V_{\Threeone}\otimes V_{\Threeone}$. Any
endomorphism $\alpha$ of $V_{\Threeone}\otimes V_{\Threeone}$ maps each space
$W_{\nu}$ to itself, and is determined by the resulting linear maps
$\alpha_{\nu}\co W_{\nu}\to W_{\nu}$. 

Where two endomorphisms $\alpha$ and $\beta$ of $\Oplus( W_{\nu}\otimes
V_{\nu})$ are composed, the corresponding restrictions to each weight space
$W_{\nu}$ compose, to give $(\alpha \circ \beta)_{\nu}=\alpha_{\nu}\circ
\beta_{\nu}$. Now the invariant of the closure of a tangle represented by an
endomorphism
$\gamma$ of  $\Oplus (W_{\nu}\otimes V_{\nu})$ is known to be $\sum
(tr(\gamma_{\nu})\times \delta_{\nu})$, where $\delta_{\nu}=J_O(V_{\nu})$ is the
quantum dimension of the module $V_{\nu}$. The difference of the invariants for
two knots represented respectively by $\gamma$ and $\gamma '$ is then given in
the same way using $\gamma-\gamma '$ in place of $\gamma$. 

The invariants for Conway and Kinoshita--Teresaka arise in this way from
endomorphisms $\gamma=\alpha\circ\beta$ and $\gamma '=\alpha '\circ\beta$, in
which $\alpha$ and $\alpha '$ represent one of the 2--tangles for Conway, and
the same tangle turned over for Kinoshita--Teresaka, while the other tangle
gives the same $\beta$ in each case. We can write $\alpha
'=R^{-1}\circ\alpha\circ R$ as module endomorphisms,  where $R$ is the
$R$--matrix for $V_{\Threeone}$. Clearly, for those
$\nu$ with $\hbox{dim }W_{\nu}=1$ we will have $\alpha '_{\nu}=\alpha_{\nu}$,
and so $\gamma '_{\nu}-\gamma_{\nu}=0$. (As noted in \cite{MC}, if this happens
for all $\nu$ then the invariant cannot distinguish any mutant pair). The final
difference of invariants will thus depend only on those
$\nu$ where the summand $V_{\nu}$ has multiplicity greater than 1. In the case
here there are just two such $\nu$ and in each case the space $W_{\nu}$ has
dimension 2. The calculation then reduces to the determination of the $2\times
2$ matrices representing $\alpha_{\nu},\alpha '_{\nu}$ and $\beta_{\nu}$.
 
\subsection{Result of the explicit calculation}
 The difference between the
values of the invariant on Conway's knot and on the  Kinoshita--Teresaka knot is 
\[
\begin{array}{rl}&s^{-80}   ({s}^{8}+1  )^{2}  ({s}^{4}+1  )^{4}
  (s+1  )^{13}  (s-1  )^{13}  ({s}^{2}-s+1
  )^{3}  ({s}^{2}+s+1  )^{3}\\ &  ({s}^{6}-{s}^{5}+{s}^{
4}-{s}^{3}+{s}^{2}-s+1  )  ({s}^{6}+{s}^{5}+{s}^{4}+{s}^{3}+{ s}^{2}+s+1  )\\
&   ({s}^{4}-{s}^{3}+{s}^{2}-s+1  )  ({s} ^{4}+{s}^{3}+{s}^{2}+s+1  ) 
({s}^{4}-{s}^{2}+1  )
  ({s}^{2}+1  )^{6}\\ &  ({s}^{46}-{s}^{44}+2\,{s}^{40}-4
\,{s}^{38}+2\,{s}^{36}+3\,{s}^{34}-4\,{s}^{32}+6\,{s}^{30}-{s}^{28}-3
\,{s}^{26}+6\,{s}^{24}\\ & -4\,{s}^{22}+4\,{s}^{20} +2\,{s}^{18}-5\,{s}^{16}
+5\,{s}^{14}-2\,{s}^{12}-2\,{s}^{10}+4\,{s}^{8}-2\,{s}^{6}+{s}^{2}-1)\\
\end{array}
\]
 up to a power of the variable $s$.

This may be rewritten to indicate more clearly the appearance of roots of unity
as the  product of $  ({s}^{46}-{s}^{44}+2\,{s}^{40}-4
\,{s}^{38}+2\,{s}^{36}+3\,{s}^{34}-4\,{s}^{32}
 +6\,{s}^{30}-{s}^{28}-3
\,{s}^{26}+6\,{s}^{24}-4\,{s}^{22}+4\,{s}^{20}+2\,{s}^{18}-5\,{s}^{16}
+5\,{s}^{14}-2\,{s}^{12}-2\,{s}^{10}+4\,{s}^{8}-2\,{s}^{6}+{s}^{2}-1
  )$ with the factors $(s^8-s^{-8})^2 (s^7-s^{-7})(s^6-s^{-6})(s^5-
s^{-5})(s^4-s^{-4})^2(s^3-s^{-3})^2(s^2-s^{-2})(s-s^{-1})^3$, and a power of
$s$.

When this is written as a power series in $h$ with $s=e^{h/2}$ the first term
becomes
$7+O(h)$ and the other factors contribute $ch^{13}+O(h^{14})$, where the
coefficient $c$ is
$c=8^2.7.6.5.4^2.3^2.2$. The coefficient of $h^{13}$ in the power series
expansion of the
 $SU(3)_q$ invariant for the 15--dimensional irreducible module is thus a
Vassiliev invariant of type at most 13  which differs on the two mutant knots.

\subsection{Some background to the calculational method}

In the following section we give details of the methods used in our
calculations. We feel it is important that others can in principle check the
calculations, as we were very  much aware in setting up our initial data just
how much scope there is for error. It can easily cause problems, for example,
if some of the data is taken from one source and some from another which has
been normalised in a slightly different way. When the goal is to show that some
polynomial arising from the calculations is non-zero any  mistake is almost
bound to result in a non-zero polynomial even if the true polynomial is zero.

In our work here we have been reassured to find that the non-zero difference
polynomial above at least has some roots which could be anticipated, since the
difference must vanish at certain roots of unity. An error in the calculations
would have been likely to give a difference which did not have these roots.

The computations were done in Maple, using its polynomial handling and linear
algebra routines. In this way we avoided the need to write explicit Pascal or C
programs for matrices and polynomials, although the computations were probably
not as fast as with a compiled program. For comparison, a Maple version of the
Hecke algebra program in \cite{MS} took roughly 50 times as long as the
compiled Pascal program to calculate the Homfly polynomial of a variety of links
when tested some time ago on the same machine.

\subsection{The quantum group $SU(3)_q$}

We start from a presentation of the quantum group $SU(3)_q$ as an algebra with
six generators, $X_1^{\pm},\,X_2^{\pm},\,H_1,\,H_2$, and a description of the
comultiplication and antipode. Let $M$ be any finite-dimensional left module
over $SU(3)_q$. The action of any one of these six generators $Y$  will
determine a linear endomorphism $Y_M$ of $M$. We  build up explicit matrices
for these endomorphisms on a selection of low-dimensional modules, using the
comultiplication to deal with the tensor product of two known modules, and the
antipode to construct the action on the linear dual of a known module. We must
eventually determine the matrices $Y_M$ for  the $15$--dimensional module 
$M=V_{\Threeone}$  above, and find the $225\times225$  $R$--matrix, $R_{MM}$ 
which represents the endomorphism of $M\otimes M$ needed for crossings.

 Knowing $Y_M$ we can find the generators $Y_{MM}$ for the module $M\otimes M$,
and thus identify the highest-weight vectors for this module. We can follow the
effect of each $2$--tangle $F$ and $G$ on the highest-weight vectors when we
know how to take account of the closure of one of the strings in forming the
$2$--tangle from the $3$--braid. To do this we need the fixed element $T$ of the
quantum group, corresponding to Turaev's `enhancement' \cite{Turaev1}, which is
used in forming the `quantum trace'.

 For the quantum groups coming from the  classical Lie algebras there is a
simple prescription for $T=\exp(h\rho)$ in terms of a linear form $\rho=\sum
\mu_i H_i$, with coefficients determined by the Cartan matrix for the Lie
algebra, \cite{Kassel}. In the case of $SU(3)_q$ we have $\rho=H_1+H_2$. The
quantum dimension of any module $M$ is the trace of the matrix $T_M$
representing the action of $T$ on $M$. More generally, the effect of closing a
string which is coloured by $M$, to convert an endomorphism of $V\otimes M$
into an endomorphism of $V$, can be realised by acting on $M$ by $T$ and then
taking the partial trace of the composite linear endomorphism of $V\otimes M$.
The element $T$ is variously written as $u^{\pm1}v$ or $u^{-1}\theta$ where $v$
is Turaev's `ribbon element' representing the full twist and $u$ is constructed
directly from the universal R--matrix, \cite{Turaev}, \cite{Kassel}.

We follow Kassel in the basic description of the quantum group from
\cite{Kassel}, chapter 17, using generators $H_1$ and $H_2$ for the Cartan
sub-algebra, but with generators $X_i^\pm$ in place of $X_i$ and $Y_i$. We  use
the notation  $K_i=\exp(hH_i/4)$, and  set $a=\exp(h/4), \, s=\exp(h/2)=a^2$
and $q=\exp(h)=s^2$, unlike Kassel. The elements satisfy the commutation
relations $[H_i,H_j]=0$, $[H_i,X_j^\pm]=\pm a_{ij}X_j^\pm$, 
$[X_i^+,X_i^-]=(K_i^2-K_i^{-2})/(s-s^{-1})$, where
$(a_{ij})=\pmatrix{2&-1\cr-1&2\cr}$ is the Cartan matrix for $SU(3)$, and
also the Serre relations of degree 3 between $X_1^\pm$ and $X_2^\pm$. 

Comultiplication is given by
\[\begin{array}{rl}\Delta(H_i)&=H_i\otimes I+I\otimes H_i,\\ (\hbox{so
}\Delta(K_i)&= K_i\otimes K_i,)\\
\Delta(X_i^\pm)&=X_i^\pm\otimes K_i+ K_i^{-1}\otimes X_i^\pm,\\
\end{array}
\] and the antipode $S$ by $S(X_i^\pm) =-s^{\pm 1}X_i^\pm$, $S(H_i)=-H_i$, $
S(K_i)=K_i^{-1}$.

The fundamental $3$--dimensional module, which we denote by $E$, has a basis in
which the quantum group generators are represented by the matrices $Y_E$ as
listed here.
\[X_1^+=\pmatrix{0&1&0\cr0&0&0\cr0&0&0\cr},\
X_2^+=\pmatrix{0&0&0\cr0&0&1\cr0&0&0\cr}\]
\[ X_1^-=\pmatrix{0&0&0\cr1&0&0\cr0&0&0\cr},\
X_2^-=\pmatrix{0&0&0\cr0&0&0\cr0&1&0\cr}\]
\[H_1=\pmatrix{1&0&0\cr0&-1&0\cr0&0&0\cr},\
H_2=\pmatrix{0&0&0\cr0&1&0\cr0&0&-1\cr}.\]

For calculations we   keep track of the elements $K_i$ rather than $H_i$,
represented by
\[K_1=\pmatrix{a&0&0\cr0&a^{-1}&0\cr0&0&1\cr},\
K_2=\pmatrix{1&0&0\cr0&a&0\cr0&0&a^{-1}\cr}\] for the module $E$.

We can then write down the elements $Y_{EE}$ for the actions of the generators
$Y$ on the module $E\otimes E$, from the comultiplication formulae. The
$R$--matrix $R_{EE}$ representing the endomorphism of $E\otimes E$ which is used
for  the crossing of two strings coloured by $E$ can be given, up to a scalar,
by the prescription
\[\begin{array}{rl}R_{EE}(e_i\otimes e_j)&=e_j\otimes e_i, \hbox{ if }i>j,\\
&=s\,e_i\otimes e_i, \hbox{ if } i=j,\\ &=e_j\otimes e_i+(s-s^{-1})e_i\otimes
e_j, \hbox{ if }i<j,
\\
\end{array}
\] for basis elements $\{e_i\}$ of $E$.

We made a quick check with Maple to confirm that the matrices $Y_{EE}$ all
commute with $R_{EE}$, as they should. It can also be checked that $R_{EE}$ has
eigenvalues $s$ with multiplicity $6$ and $-s^{-1}$ with multiplicity $3$, and
satisfies the equation $R-R^{-1}=(s-s^{-1})\hbox{Id}$.

The linear dual $M^*$ of a module $M$ becomes a module when the action of a
generator $Y$ on $f\in M^*$ is defined by $<Y_{M^*}f,v>=<f,S(Y_M)v>$, for $v\in
M$. For the dual module $F=E^*$ we then have matrices for $Y_F$, relative to
the dual basis, as follows.

\[X_1^+=\pmatrix{0&0&0\cr-s&0&0\cr0&0&0\cr},\
X_2^+=\pmatrix{0&0&0\cr0&0&0\cr0&-s&0\cr}\]
\[ X_1^-=\pmatrix{0&-s^{-1}&0\cr0&0&0\cr0&0&0\cr},\
X_2^-=\pmatrix{0&0&0\cr0&0&-s^{-1}\cr0&0&0\cr}\]
\[K_1=\pmatrix{a^{-1}&0&0\cr0&a&0\cr0&0&1\cr},\
K_2=\pmatrix{1&0&0\cr0&a^{-1}&0\cr0&0&a\cr}.\]

The most reliable way to work out the $R$--matrices $R_{EF}, R_{FE}$ and
$R_{FF}$ is to combine $R_{EE}$ with module homomorphisms $\hbox{cup}_{EF}$, $
\hbox{cup}_{FE}$, $\hbox{cap}_{EF}$ and $\hbox{cap}_{FE}$ between the modules
$E\otimes F$, $F\otimes E$ and the trivial 1--dimensional module, $I$, on which
$X_i^\pm$ acts as zero and $K_i$ as the identity. For example, to represent a
homomorphism from $I$ to $E\otimes F$ the matrix for $\hbox{cup}_{EF}$ must
satisfy $Y_{EF}\,\hbox{cup}_{EF}=\hbox{cup}_{EF}\,Y_{I}$, which identifies
$\hbox{cup}_{EF}$ as a common eigenvector of the matrices $Y_{EF}$, with
eigenvalue $0$ or $1$ depending on $Y$. The matrices are determined up to a
scalar by such considerations; when one has been chosen the scalar for the
others is dictated by diagrammatic considerations. They  are quite easy
to write down theoretically, although to be careful about compatibility and
possible miscopying it is as well to get Maple to find them in this way for
itself.
Once these matrices have been found they can be combined with the matrix
$R_{EE}^{-1}$ to construct the $R$--matrices $R_{EF},R_{FE},R_{FF}$, using the
diagram shown in figure 3, for example, to determine $R_{EF}$. This gives 
\[ R_{EF}=1_F\otimes 1_E\otimes  \hbox{cap}_{EF}\circ 
 1_F\otimes R_{EE}^{-1}\otimes 1_F\circ \hbox{cup}_{FE}\otimes
1_E\otimes 1_F .\]

\begin{figure}[htbp]
\cl{%
\SetLabels
\scriptsize
\E(.05*.8) $F$\\
\E(.23*.8) $E$\\
\E(.06*.23) $E$\\
\E(.23*.23) $F$\\
\E(.6*.95) $F$\\
\E(.7*.95) $E$\\
\E(.56*.69) $F$\\
\E(.74*.69) $E$\\
\E(.83*.69) $E$\\
\E(.98*.69) $F$\\
\E(.58*.42) $F$\\
\E(.72*.42) $E$\\
\E(.83*.42) $E$\\
\E(.98*.42) $F$\\
\E(.84*.05) $E$\\
\E(.95*.05) $F$\\
\E(.4*.5) {\Large $=$}\\
\endSetLabels
%\ShowGrid
\AffixLabels{\BoxedEPSF{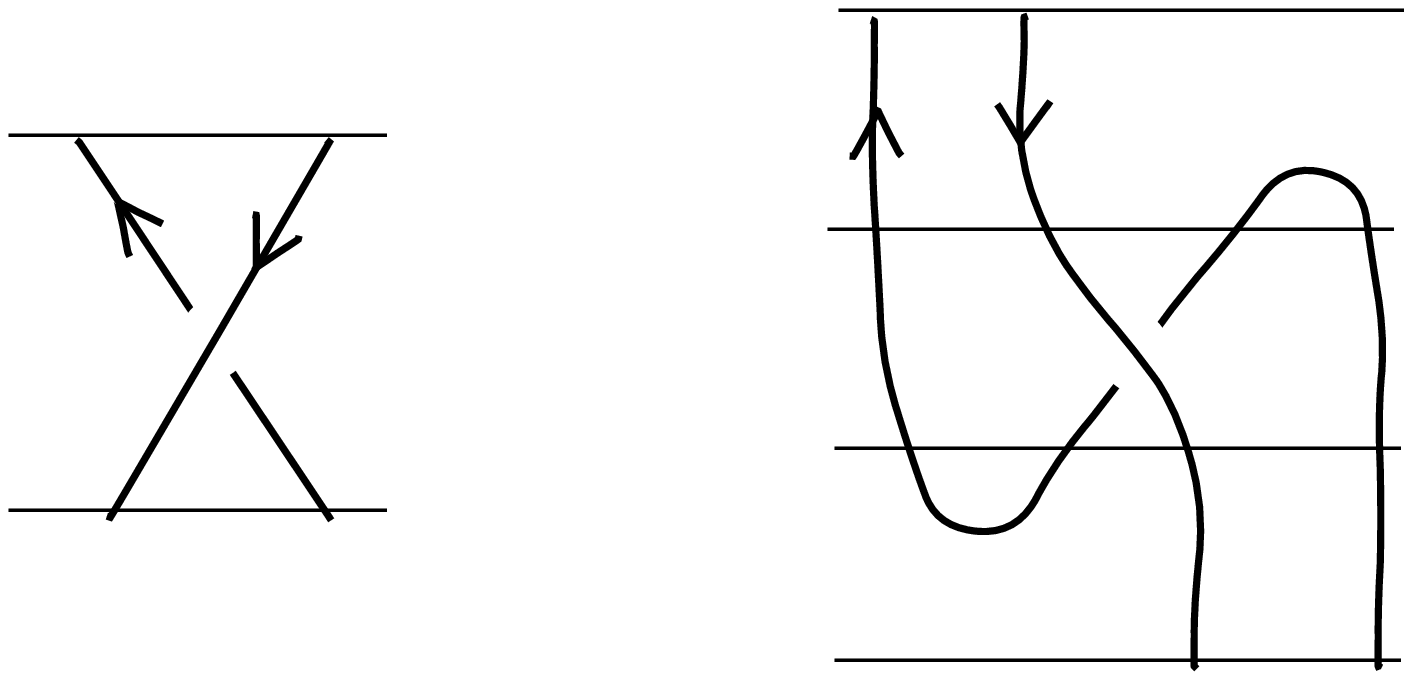 scaled  500}}}
\cl{\small Figure 3}
\end{figure}

 The module structure
of 
$M=V_{\Threeone}$ can  be found by identifying $M$ as a
$15$--dimensional submodule of $E\otimes E\otimes F$. We know that there will be
a direct sum decomposition of $E\otimes E\otimes F$ as $M\oplus N$, and indeed
that $N$ will decompose further into the sum of two copies of a $3$--dimensional
module isomorphic to $E$ and one $6$--dimensional module with Young diagram
\Twotwo.
 The full twist element on the three strings coloured by $E, E$ and $F$ acts by
a scalar on each of the irreducible submodules of $E\otimes E\otimes F$. It can
be expressed  as a $27\times27$ matrix in terms of  the $R$--matrices above.
Maple can then produce a basis for each of the eigenspaces, one of dimension 15
and the other two each of dimension 6. Write $P$ and $Q$ for the $27\times15$
and
$27\times12$ matrices whose columns are made of these basis vectors. Then $P$
and $Q$ give bases for $M$ and $N$ respectively. The partitioned matrix $(P|Q)$
is invertible. When its inverse, found by Maple, is written in the form
$\ds\left(\ds{R\over S}\right)$ we have a $15\times27$ matrix $R$ which
satisfies $RP=I_{15}$ and $RQ=0$. Regard $P$ as the matrix representing the
inclusion of the module $M$ into $E\otimes E\otimes F$. Then  $R$ is the
matrix, in the same basis, of the projection from $E\otimes E\otimes F$ to $M$.
The module generators $Y_M$   satisfy  $Y_M=R\,Y_{EEF}\,P$, giving the explicit
action of the quantum group on $M$. 

We use the injection and projection further to find the $15^2\times15^2$
$R$--matrix $R_{MM}$.  First include $M\otimes M$ in $(E\otimes E\otimes
F)\otimes(E\otimes E\otimes F)$, then construct the $R$--matrix for $E\otimes
E\otimes F$ from the crossing of three  strings each coloured with $E$ or $F$
over three others using the various matrices $R_{EF}$ from above, and finally
project to $M\otimes M$. 

The calculations can be completed in principle from here.  Represent the
$3$--braid in the $2$--tangle $F$ by an endomorphism of $M\otimes M\otimes M$,
using $R_{MM}$ and its inverse. Then use $T_M$ and the partial trace to close
off one string, hence giving the endomorphism $F_{MM}$ of $M\otimes M$
determined by $F$. A similar calculation gives the endomorphism  $G_{MM}$.  The
invariant for one of the knots is given by the trace of $T_{MM}F_{MM}G_{MM}$.
The other is given by replacing $G_{MM}$ with the conjugate
$R_{MM}^{-1}G_{MM}R_{MM}$. Some calculation can be avoided by using
$G_{MM}-R_{MM}^{-1}G_{MM}R_{MM}$ in place of $G_{MM}$, to get the difference of
the invariants directly. 

A considerable shortcut can be made at this point by concentrating  on the
effect of $F_{MM}$ and $G_{MM}$ on certain highest weight vectors in $M\otimes
M$, rather than considering the whole of the module. A {\sl highest weight
vector}
$v$ of a module $V$ is a common eigenvector of $H_1$ and $H_2$ (or equally
$K_1$ and $K_2$) which satisfies $X_1^+(v)=X_2^+(v)=0$. The submodule of $V$
generated by a highest weight vector is irreducible. Its isomorphism type is
determined by the eigenvalues of $H_1$ and $H_2$, which are non-negative
integers. It follows easily from the relations in the quantum group that any
module homomorphism $f\co V\to W$ carries highest weight vectors to highest weight
vectors of the same type.

 Calculation in Maple determines the linear subspace of $M\otimes M$ which is
the common null-space of $X_1^+$ and $X_2^+$. This turns out to have dimension
$10$, spanned by two highest weight vectors of type $(3,1)$, two of type
$(1,2)$ and six further highest weight vectors each of a different type. Then
the endomorphism $F$  restricts to a linear endomorphism $F_\nu$ of the space
of highest weight vectors of type $\nu$, for each $\nu$. We remarked earlier
that weight spaces of dimension 1 will not contribute to the difference of the
invariants on two mutant knots, so we need only calculate the maps $F_\nu$ and
$G_\nu$ for the two 2--dimensional weight spaces $\nu=(3,1)$ and $\nu=(1,2)$. We
thus choose two spanning vectors for one of these spaces and follow each of
these through the 2--tangle $F$, taking the tensor product with $M$ and mapping
to $M\otimes M\otimes M$ as above (using repeated operations of the
$225\times225$ $R$--matrix on a vector of length $225\times15$) before applying
the matrix $T_M$ and taking a partial trace to finish in $M\otimes M$. Since
the result in each case must be a linear combination of the two chosen weight
vectors it is not difficult to find the exact combination. This  determines a
$2\times2$ matrix representing $F_\nu$ for the weight space of type $\nu$.
Similar calculations for the other weight space and for $G$, along with a quick
calculation of the $2\times2$ matrix representing $R_{MM}$ on each weight type
gives enough to find the contribution of each of these weight types to the
difference. The final difference comes from multiplying the trace of the
$2\times2$ difference matrix for each type $\nu$  by the quantum dimension of
the irreducible module of type $\nu$ for each of the two types and then adding
the results. 

Up to the same power of $s$ in each case the contribution from the weight space
of type $(3,1)$ was found to be
\[\begin{array}{rl} t_{31}=&
   ({s}^{8}+1  )^{2}  ({s}^{2}+1  )^{4}
  ({s}^{4}+1  )^{3}  (s+1  )^{13}  (s-1  ) ^{13}{s}^{6}  ({s}^{2}-s+1  ) 
({s}^{2}+s+1  )\\ &
  ({s}^{4}-{s}^{3}+{s}^{2}-s+1  )  ({s}^{4}+{s}^{3}+{s}^{ 2}+s+1  ) \\ &
({s}^{6}-{s}^{5}+{s}^{4}-{s}^{3}+{s}^{2}-s+1
  )  ({s}^{6}+{s}^{5}+{s}^{4}+{s}^{3}+{s}^{2}+s+1  )\\
&(2\,{s}^{20}+{s}^{18}+{s}^{14}-{s}^{12}+2\,{s}^{8}-{s}^{6}-1
  ) \\ &
  ({s}^{22}-{s}^{20}+{s}^{16}-2\,{s}^{14}+3\,{s}^{12}+2\,{s}^{10}-
{s}^{8}+2\,{s}^{6}+2  )\\
=&(2\,{s}^{20}+{s}^{18}+{s}^{14}-{s}^{12}+2\,{s}^{8}-{s}^{6}-1
  ) \\ &
  ({s}^{22}-{s}^{20}+{s}^{16}-2\,{s}^{14}+3\,{s}^{12}+2\,{s}^{10}-
{s}^{8}+2\,{s}^{6}+2  )\\ &\times(s^8-s^{-8})^2 (s^7-s^{-7})(s^5-
s^{-5})(s^4-s^{-4})\\ &(s^3-s^{-3})(s^2-s^{-2})(s-s^{-1})^6s^{49},\\
\end{array}
\] and the contribution from type $(1,2)$ to be

\[
\begin{array}{rl} t_{12}= & ({s}^{6}-{s }^{5}+{s}^{4}-{s}^{3}+{s}^{2}-s+1 
)^{2}  ({s}^{6}+{s}^{5}+{s }^{4}+{s}^{3}+{s}^{2}+s+1  )^{2}\\ &
  ({s}^{4}-{s}^{2}+1  )
  ({s}^{8}+1  )^{2}   ({s}^{4}+1  )^{5}  ({s}^{2}+1  )^{8} \\ & ({s}^{2}+s+1 
)  ({s}^{2}-s+1  ) (s-1  )^{14}
  (s+1  )^{14} ({s}^{10}-{s}^{8}+{s}^{4}-{s}^{2}+1  )\\ &
({s}^{18}-{s}^{16}-{s}^{14}+2\,{s}^{
12}-2\,{s}^{10}+2\,{s}^{6}-2\,{s}^{4}-{s}^{2}+1  )\\ =&
({s}^{18}-{s}^{16}-{s}^{14}+2\,{s}^{
12}-2\,{s}^{10}+2\,{s}^{6}-2\,{s}^{4}-{s}^{2}+1  )\\ &
({s}^{10}-{s}^{8}+{s}^{4}-{s}^{2}+1  )\\ &\times (s^8-s^{-8})^2
(s^7-s^{-7})^2(s^6-s^{-6})(s^4-s^{-4})^3\\
&\times(s^2-s^{-2})^2(s-s^{-1})^4s^{56}.\\
\end{array}
\] 

The quantum dimension for the irreducible module of type $(3,1)$, which has
Young diagram \Fourone, is a product of quantum integers
$[6][4]=(s^6-s^{-6})(s^4-s^{-4})/(s-s^{-1})^2$. For the module of type $(1,2)$,
with Young diagram \Threetwo, it is
$[5][3]=(s^5-s^{-5})(s^3-s^{-3})/(s-s^{-1})^2$. 

The difference between the $SU(3)_q$ invariants with the module $V_{\Threeone}$
for the Conway and Kinoshita--Teresaka knots is then given, up to a power of
$s=e^{h/2}$, by $[5][3]t_{12}+[6][4]t_{31}$. This yields  the polynomial quoted
earlier.

\section{The Kuperberg skein for mutants} Let $K$ and $K'$ be the mutants shown
schematically in figure 1.  As $K$ and $K'$ are knots,  precisely one of $F$ or
$G$
 must induce the identity permutation on the endpoints
 by following the strings through the tangle, while the other induces
the transposition. We will
consider these two cases separately.

In \cite{K} Kuperberg gives  a skein-theoretic method for handling the
$SU(3)_q$ invariant of a link when coloured by the fundamental module, which he
denotes by $<>_{A_2}$. Knot diagrams are extended to allow 3--valent oriented
graphs in which any vertex is either a sink or a source. Crossings can be
replaced locally in this skein  by  a  linear combination of planar graphs, and
any planar circles, 2--gons or 4--gons can be replaced by linear combinations of
simpler pieces.

In using skein-based calculations it is helpful when dealing, for example, with
satellites to regard the pattern as a diagram in an annulus, and note that it
can be replaced by any equivalent linear combination of diagrams in the skein
of the annulus. Thus we should consider the Kuperberg skein of the annulus,
namely linear combinations of admissibly oriented 3--valent graph diagrams
subject to local relations as  before. A similar definition can be made for the
skein of other surfaces. Notice that the relations ensure that the skein is
spanned by oriented graphs lying entirely in the surface, without simple closed
curves, 2--gons or 4--gons which bound discs in the surface. 

In the case of the annulus this shows that the skein is spanned by unions of
oriented simple closed curves parallel to the boundary of the annulus, with
orientations in either direction. 

When a mutant knot $K$ is made up from two $2$--tangles $F$ and $G$ as above
then one of $F$ and $G$, let us suppose $G$, must be a pure tangle, in the
sense that the arcs of $G$ connect the entry point at top left with the exit at
bottom left, and top right with bottom right. Then $K$ can be viewed as made
from the diagram in the disc $P$ with two holes, shown in figure 4, by embedding
the planar surface $P$ so that  the two `ears' are tied around the arcs of $G$.
Turning the diagram in $P$ over along the axis indicated before embedding
it in the same way, and reversing all string orientations, will give one of the
mutants
$K'$ of
$K$. Any satellites of $K$ and $K'$ are related in a similar way, for we can
view a satellite of $K$ as constructed by decorating the diagram in $P$ with
the required pattern, and then tying the ears of $P$ around $G$ as before. The
corresponding satellite of $K'$ is given by turning $P$ over, with the
decorated diagram, reversing all strings, and then using the same embedding of
$P$. 

\begin{figure}[htbp]
\cl{%\small
\SetLabels
\E(.54*.51) $F$\\
\endSetLabels
%\ShowGrid
\AffixLabels{$P\ =\quad${\BoxedEPSF{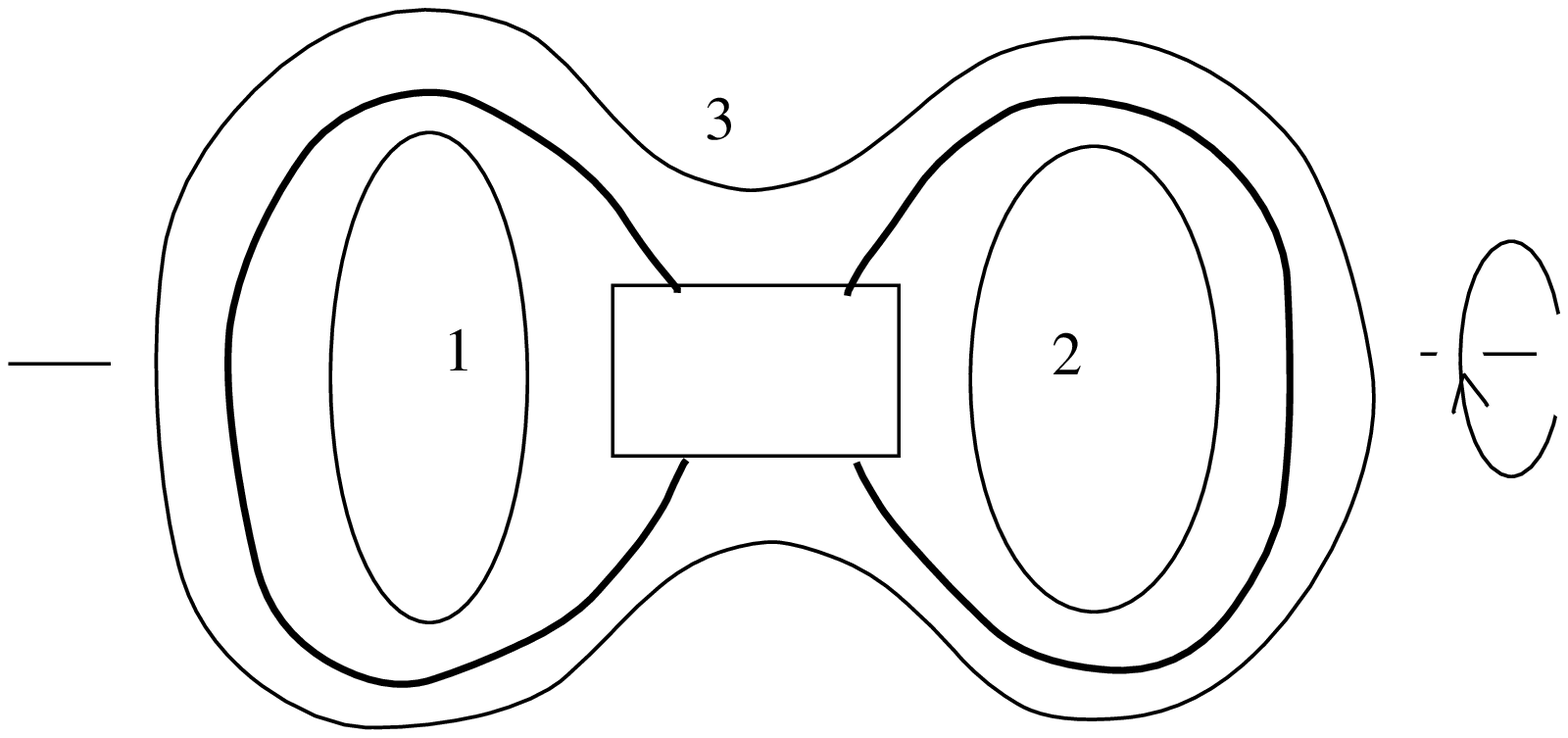 scaled  500}}}}
\cl{\small Figure 4}
\end{figure}

If we could show that the Kuperberg skein of $P$ is spanned by elements which
are invariant under turning over and reversing orientation then we could deduce
that satellites of mutants such as $K$ and $K'$ would have the same $SU(3)_q$
invariants, by considering the decorated diagram in this skein.  A proof for
all mutants would need a similar analysis for the skein of the once-punctured
torus, to deal with one of the other mutation operations, and the third case
would then follow, using a similar argument to \cite{MT}, where the truth of
the corresponding results in the Kauffman bracket skein showed that satellites
of mutants have the same $SU(2)_q$ invariants.

We shall now describe a basis for the Kuperberg skein of $P$, which has
some nice symmetry properties, but not enough to give the invariance above.
Indeed a diagram coming from a $3$--fold parallel with one reversed string will
give a linear combination of basis elements in the skein in which all but at
most one pair are invariant. (Diagrams from $2$--fold parallels of any
orientation determine elements of the invariant subspace.)

\begin{thm}\sl The Kuperberg skein of a disc with two holes  has a basis of
diagrams consisting of the union of simple closed curves parallel to each
boundary component and a trivalent graph with a $2$--gon nearest to each of the
three boundary components and $6$--gons elsewhere.\end{thm}
 
\begin{proof} Use the skein relations to write any diagram as a
linear combination of admissibly oriented trivalent graphs in the surface. We
can assume that there are no simple closed curves or $2$--gons or $4$--gons with 
null-homotopic boundary. There may be a number of simple closed curves parallel
to each of the boundary components. The remaining graph must be connected, otherwise 
one of its components lies in an annulus inside the surface, and
can be reduced further to a linear combination of unions of parallel simple
closed curves. Consider the graph  as lying in $S^2$, by filling in the three
boundary components of the surface. It dissects $S^2$ into a number of
$n$--gons, with $n$ even, and $n\ge6$ except possibly for the three $n$--gons
containing the added discs.
 Now calculate the Euler characteristic of the resulting sphere $S$ from the
dissection by the graph. As vertices are trivalent and each edge now bounds two
faces, we can count the Euler characteristic as a sum over the $n$--gons, in
which each vertex  contributes $ 1/3$ and  each edge
$-1/2$. Therefore each $n$--gon will contribute $1 - n/6$,  so the only positive
contribution to $\chi(S)$ can come from $2$--gons or $4$--gons.  These can only
arise from the original three boundary components, where the maximum possible
total positive contribution is $2$ when each boundary component gives a
$2$--gon. Since the total must be $2$ and the only other contributions are
negative or zero,
 we must have three $2$--gons forming the original boundary components and
$6$--gons elsewhere.
\end{proof}

If we start with a $3$--parallel of a tangle $F$ inside the planar surface $P$,
with two strands in one direction and one in the other, and write it in the
Kuperberg skein we will get a linear combination of graphs as above, each
having at most $3$ strings around each `ear'. Some of these will be the union
of some simple closed curves around the punctures and trivalent graphs. In
figure 5 we show  one such trivalent graph which fails to be symmetric
under the order 2 operation of turning the surface over (and reversing edge
orientations). 

\begin{figure}[htbp]
\cl{\Skeinelement}
\cl{\small Figure 5}
\end{figure}

Note however that this graph is symmetric under the operation of order 3 in
which the three boundary components are cycled. This is 
a general feature of the connected trivalent graphs which arise in our
construction, as appears from the following description, where we replace $P$
by a $3$--punctured sphere. 

 We call a trivalent graph in the
$3$--punctured sphere
 {\it admissible} if it is oriented so that each vertex is either a sink or a
source, and every region not containing a puncture is a hexagon. 

\begin{thm}\sl Every admissible graph in the $3$--punctured sphere is
symmetric, up to isotopy avoiding the punctures, under a rotation which cycles
the punctures. It can be constructed from the hexagonal tesselation of the
plane by choosing an equilateral triangle lattice whose vertices lie at the
centres of some of the hexagons and factoring out the translations of the
lattice and the  rotations of order $3$ which preserve the lattice.
\end{thm}

\begin{proof} Let $\Gamma$ be the admissible graph. By our Euler characteristic
calculations we know that each puncture is contained in a $2$--gon. There is a
$3$--fold branched cover of $S^2$ by the torus $T^2$ with three branch points,
each  cyclic of order $3$. The inverse image of $\Gamma$ in $T^2$ then consists
of hexagonal regions, with three distinguished regions containing the branch
points. This inverse image is invariant under the deck transformation of order
$3$ which leaves each distinguished region invariant. The further inverse image
under the regular covering of $T^2$ by the plane is a tesselation of the plane
by hexagons, and the inverse image of the centre of one of the distinguished
regions determines a lattice in the plane. We want to show that this is an
equilateral triangle lattice, when the hexagonal tesselation is drawn in the
usual way. We need only lift the deck transformation to a transformation of the
plane keeping the tesselation invariant and fixing one of the lattice points to
see that it must lift to a rotation of the tesselation about the centre of a
distinguished hexagon. Since the lattice is invariant under this transformation
it follows that the lattice must be equilateral. The inverse image of each of
the other two branch points will also form an equilateral lattice, invariant
under the first rotation, and so their vertices lie in the centres of the
triangles; by construction they also lie in the middle of hexagons. Although
the equilateral lattice need not lie symmetrically with respect to reflections
of the tesselation, as in the example shown below,  it does follow that the
rotation which permutes the three lattices will also preserve the tesselation.
This rotation induces the symmetry of the sphere which cycles the branch points
and preserves $\Gamma$.
\end{proof}

\begin{figure}[htbp]
\cl{\Hexnet}
\vskip .1 in
\cl{\small Figure 6}
\end{figure}

Figure 6 shows such an equilateral triangle lattice superimposed on a hexagon
tesselation. 
The resulting graph in the 3--punctured sphere, whose
fundamental domain is indicated, is the graph shown in figure 5 as a
non-symmetric skein element in the disk with  two holes. The labelling of the
puncture points as 1, 2 and 3 corresponds to that of the boundary components.
The 3--fold symmetry of the graph in the surface when the boundary components
are cycled is  evident from this viewpoint.

\medskip
The Kuperberg skein of the punctured torus does not appear to have such a
simple basis. The region around the puncture may be a $2$--gon or a
$4$--gon, giving the  following possible combinations: (i) a $2$--gon, two
$8$--gons and $6$--gons elsewhere,  (ii) a $2$--gon, one $10$--gon and $6$--gons
elsewhere, (iii) a $4$--gon, one $8$--gon and $6$--gons elsewhere, (iv) $6$--gons
only. We did not try to analyse the configurations further, in view of the
results of our quantum calculations.

\newpage
%%%%%%%%%%%%%%%%%%%%%  End of main body of article

\Addresses\recd

\end{document}